\documentclass[a4paper,11pt]{amsart}
\usepackage{graphicx}
\usepackage{mathptmx}
\usepackage{amsmath}
\usepackage{amssymb}
\usepackage{enumitem}
\usepackage{xcolor}
\usepackage{xparse}
\NewDocumentCommand{\eulerian}{omm}
 {%
  \genfrac<>{0pt}{}{#2}{#3}%
  \IfValueT{#1}{_{\!#1}}%
 }

 \newmuskip\pFqmuskip

\newcommand*\pFq[6][8]{%
  \begingroup 
  \pFqmuskip=#1mu\relax
  \mathchardef\normalcomma=\mathcode`,
  \mathcode`\,=\string"8000
  \begingroup\lccode`\~=`\,
  \lowercase{\endgroup\let~}\pFqcomma
  {}_{#2}F_{#3}{\left(\genfrac..{0pt}{}{#4}{#5}\bigg|#6\right)}%
  \endgroup
}
\newcommand{\pFqcomma}{{\normalcomma}\mskip\pFqmuskip}

\begin{document}

\title[Normal ordering of degenerate integral powers of number operator]{Normal ordering of degenerate integral powers of number operator and its applications}

\author{Taekyun  Kim}
\address{Department of Mathematics, Kwangwoon University, Seoul 139-701, Republic of Korea}
\email{tkkim@kw.ac.kr}

\author{Dae San Kim}
\address{Department of Mathematics, Sogang University, Seoul 121-742, Republic of Korea}
\email{dskim@sogang.ac.kr}

\author{Hye Kyung Kim}
\address{Department of Mathematics Education, Daegu Catholic University, Gyeonsang, Republic of Korea}
\email{hkkim@cu.ac.kr}


\maketitle

\begin{abstract}
The normal ordering of an integral power of the number operator in terms of boson operators is expressed with the help of the Stirling numbers of the second kind. As a `degenerate version' of this, we consider the normal ordering of a degenerate integral power of the number operator in terms of boson operators, which is represented by means of the degenerate Stirling numbers of the second kind. As an application of this normal ordering, we derive two equations defining the degenerate Stirling numbers of the second kind and a Dobinski-like formula for the degenerate Bell polynomials.
\end{abstract}

\section{Introduction}

The Stirling number of the second $S_{2}(n,k)$ is the number of ways to partition a set of $n$ objects into $k$ nonempty subsets. The Stirling numbers of the second kind have been extensively studied and repeatedly and independently rediscovered during their long history. The Stirling numbers of the second kind appear in many different contexts and have numerous applications, for example to enumerative combinatorics and quantum mechanics. They are given either by \eqref{5} or by \eqref{7}. The study of degenerate versions of some special numbers and polynomials began with Carlitz's paper in [2], where the degenerate Bernoulli and Euler numbers were investigated. It is remarkable that in recent years quite a few degenerate versions of special numbers and polynomials have been explored with diverse tools and yielded many interesting results (see [4-6] and the references therein). It turns out that the degenerate Stirling numbers of the second play an important role in this exploration for degenerate versions of many special numbers and polynomials. The normal ordering of an integral power of the number operator $a^{\dagger}a$ in terms of boson operators $a$ and $a^{\dagger}$ can be written in the form
\begin{equation*}
(a^{\dagger}a)^{k}=\sum_{l=0}^{k}S_{2}(k,l)(a^{\dagger})^{l}a^{l}.
\end{equation*}
In addition, the normal ordering of the degenerate $k$th power of the number operator $a^{\dagger}a$, namely $(a^{\dagger}a)_{k,\lambda}$, in terms of boson operators $a, a^{\dagger}$ can be written in the form
\begin{equation}
(a^{\dagger}a)_{k,\lambda}=\sum_{l=0}^{k}S_{2,\lambda}(k,l)(a^{\dagger})^{l}a^{l},\label{1}
\end{equation}
where the generalized falling factorials $(x)_{n,\lambda}$ are given by \eqref{3} and the degenerate Stirling numbers $S_{2,\lambda}(k,l)$ by \eqref{4} and \eqref{6}.\par
The aim of this paper is to use the normal ordering in \eqref{1} in order to derive two equations defining the degenerate Stirling numbers of the second kind (see \eqref{4}, \eqref{6}) and a Dobinski-like formula for the degenerate Bell numbers (see \eqref{37}). In more detail, our main results are as follows. Firstly, by applying the degenerate $k$th power $(a^{\dagger}a)_{k,\lambda}$ to any number states $|m\rangle,\ m=0,1,2,\dots$ and using \eqref{1} we obtain the equation \eqref{4}. Secondly, by comparing one expression of $\langle z|e_{\lambda}^{a^{\dagger}a}(t)|z \rangle$ obtained by using \eqref{1} and another one of it obtained by solving a differential equation we get the equation \eqref{6}.
Thirdly, we obtain a Dobinski-like fromula for the degenerate Bell numbers $\phi_{k,\lambda}$ by showing that $\langle z|(a^{\dagger}a)_{k,\lambda}|z \rangle$ is equal to $\phi_{k,\lambda}(|z|^{2})$ by using \eqref{1} and that it is also equal to some other expression coming from the representation of coherent  state in terms of number states. For the rest of this section, we recall what are needed throughout this paper.\par
For any $\lambda\in\mathbb{R}$, the degenerate exponential functions are defined by
\begin{equation}
e_{\lambda}^{x}(t)=\sum_{k=0}^{\infty}\frac{(x)_{n,\lambda}}{k!}t^{k},\quad (\mathrm{see}\ [2,4,5,6]), \label{2}
\end{equation}
where the generalized falling factorials $(x)_{n,\lambda}$ are defined by
\begin{equation}
(x)_{0,\lambda}=1,\quad (x)_{n,\lambda}=x(x-\lambda)\cdots\big(x-(n-1)\lambda\big),\quad (n\ge 1).\label{3}
\end{equation}
When $x=1$, we let $\displaystyle e_{\lambda}(t)=e_{\lambda}^{1}(t)=\sum_{k=0}^{\infty}\frac{(1)_{k,\lambda}}{k!}t^{k}.\displaystyle$ The degenerate Stirling numbers of the second kind are defined by
\begin{equation}
(x)_{n,\lambda}=\sum_{k=0}^{n}S_{2,\lambda}(n,k)(x)_{k},\quad (n\ge 0),\quad (\mathrm{see}\ [4]),\label{4}
\end{equation}
where $(x)_{0}=1,\ (x)_{n}=x(x-1)\cdots(x-n+1),\ (n\ge 1)$. \\
Note that $\displaystyle\lim_{\lambda\rightarrow 0}S_{2,\lambda}(n,k)=S_{2}(n,k)\displaystyle$ are the ordinary Stirling numbers of the second kind given by
\begin{equation}
x^{n}=\sum_{k=0}^{n}S_{2}(n,k)(x)_{k},\quad (n\ge 0),\quad (\mathrm{see}\ [1,4]). \label{5}
\end{equation}
From \eqref{4}, we note that
\begin{equation}
\frac{1}{k!}\big(e_{\lambda}(t)-1\big)^{k}=\sum_{n=k}^{\infty}S_{2,\lambda}(n,k)\frac{t^{n}}{n!},\quad (k\ge 0),\quad (\mathrm{see}\ [4]). \label{6}
\end{equation}
By letting $\lambda \rightarrow 0$ in \eqref{6}, we see that the Stirling numbers of the second kind are also given by
\begin{equation}
\frac{1}{k!}\big(e^{t}-1\big)^{k}=\sum_{n=k}^{\infty}S_{2}(n,k)\frac{t^{n}}{n!},\quad (k\ge 0),\quad (\mathrm{see}\ [4]). \label{7}
\end{equation}
In [6], the degenerate Bell polynomials are defined by
\begin{equation}
e^{x(e_{\lambda}(t)-1)}=\sum_{n=0}^{\infty}\phi_{n,\lambda}(x)\frac{t^{n}}{n!}.\label{8}	
\end{equation}
Thus, by \eqref{8}, we get
\begin{displaymath}
	\phi_{n,\lambda}(x)=\sum_{k=0}^{n}S_{2,\lambda}(n,k)x^{k},\quad (n\ge 0),\quad (\mathrm{see}\ [5,6]).
\end{displaymath}
When $x=1$, $\phi_{n,\lambda}=\phi_{n,\lambda}(1)$ are called the degenerate Bell numbers. \par
From \eqref{4}, we note that
\begin{equation}
S_{2,\lambda}(n+1,k)=S_{2,\lambda}(n,k-1)+(k-n\lambda)S_{2,\lambda}(n,k),\quad (\mathrm{See}\ [4]),\label{9}
\end{equation}
where $n,k\in\mathbb{N}$ with $n\ge k$. \par
In this paper, we pay attention to some properties with the boson operators $a, a^{\dagger}$ that satisfy
\begin{displaymath}
[a,a^{\dagger}]=aa^{\dagger}-a^{\dagger}a=1,\quad (\mathrm{see}\ [3,7]).
\end{displaymath}
The normal ordering of an integral power of the number operator $a^{\dagger}a$ in terms of boson operators $a$ and $a^{+}$ can be written in the form
\begin{equation}
(a^{\dagger}a)^{k}=\sum_{l=0}^{k}S_{2}(k,l)(a^{\dagger})^{l}a^{l},\quad (\mathrm{see}\ [3,7]).\label{10}
\end{equation}
The number states $|m\rangle,\ m=0,1,2,\dots$, are defined as
\begin{equation}
a|m\rangle=\sqrt{m}|m-1\rangle,\quad a^{\dagger}|m\rangle=\sqrt{m+1}|m+1\rangle.\label{11}
\end{equation}
By \eqref{11}, we get $a^{\dagger}a|m\rangle=m|m\rangle$. The coherent states $|z\rangle$, where $z$ is a complex number, satisfy $a|z\rangle=z|z\rangle,\quad \langle z|z\rangle=1$.  To show a connection to coherent states, we recall that the harmonic oscillator has Hamiltonian $H=a^{\dagger}a$ (neglecting the zero point energy) and the usual eigenstates $|n\rangle$ (for $n\in\mathbb{N}$) satisfying $H|n\rangle=n|n\rangle$ and $\langle m|n\rangle=\delta_{m,n}$, where $\delta_{m,n}$ is Kronecker's symbol. \par
In this paper, we show that the combinatorial properties of the degenerate Bell numbers and polynomials are derived from the algebraic properties of the boson operators.

\section{Normal ordering of degenerate integral powers of the number operator and its applications}
First, we recall the definition of coherent states $a|z\rangle=z|z\rangle$, equivalently $\langle z|a^{\dagger}=\langle z|\overline{z}$, where $z\in\mathbb{C}$. For the coherent states $|z\rangle$, we write
\begin{equation}
|z\rangle	=\sum_{n=0}^{\infty}A_{n}|n\rangle.\label{12}
\end{equation}
Then, by \eqref{12}, we get
\begin{equation}
a|z\rangle=\sum_{n=0}^{\infty}A_{n}a|n\rangle=\sum_{n=1}^{\infty}A_{n}\sqrt{n}|n-1\rangle,	\label{13}
\end{equation}
and
\begin{equation}
a|z\rangle=z|z\rangle=z\sum_{n=0}^{\infty}A_{n}|n\rangle=\sum_{n=1}^{\infty}zA_{n-1}|n-1\rangle.\label{14}
\end{equation}
From \eqref{13} and \eqref{14}, we have
\begin{align}
A_{n}&=\frac{z}{\sqrt{n}}A_{n-1}=\frac{z}{\sqrt{n}}\frac{z}{\sqrt{n-1}}A_{n-2}=\cdots
=\frac{z^{n}}{\sqrt{n!}}A_{0}.\label{15}
\end{align}
So, by \eqref{12} and \eqref{15}, we get
\begin{equation}
|z\rangle=A_{0}\sum_{n=0}^{\infty}\frac{z^{n}}{\sqrt{n!}}|n\rangle.\label{16}	
\end{equation}
By the property of coherent state $|z\rangle$, we get
\begin{align}
1&=\langle z|z\rangle=\overline{A}_{0}\sum_{m=0}^{\infty}\frac{\overline{z}^{m}}{\sqrt{m!}} A_{0}\sum_{n=0}^{\infty}\frac{z^{n}}{\sqrt{n!}}\langle m|n\rangle\label{17} \\
&=|A_{0}|^{2}\sum_{n=0}^{\infty}\frac{|z|^{2n}}{n!}=|A_{0}|^{2}e^{|z|^{2}}.\nonumber	
\end{align}
Thus, by \eqref{17}, we get
\begin{displaymath}
A_{0}=e^{-\frac{1}{2}|z|^{2}}\cdot e^{i(\textnormal{artitrary phase})}.
\end{displaymath}
Discarding the phase factors, by \eqref{16}, we get
\begin{equation}
|z\rangle=e^{-\frac{|z|^{2}}{2}}\sum_{n=0}^{\infty}\frac{z^{n}}{\sqrt{n!}}|n\rangle.\label{18}
\end{equation}
For $x,y\in\mathbb{C}$, we have
\begin{align*}
	\langle x|y\rangle &=e^{-\frac{|x|^{2}}{2}}\sum_{m=0}^{\infty}\frac{(\overline{x})^{m}}{\sqrt{m!}} e^{-\frac{|y|^{2}}{2}}\sum_{n=0}^{\infty}\frac{y^{n}}{\sqrt{n}!}\langle m|n\rangle \\
	&=e^{-\frac{|x|^{2}}{2}-\frac{|y|^{2}}{2}}\sum_{n=0}^{\infty}\frac{(\overline{x}y)^{n}}{n!}=e^{-\frac{1}{2}(|x|^{2}+|y|^{2})+\overline{x}y}.
\end{align*}
We recall that the standard bosonic commutation relations $[a,a^{\dagger}]=aa^{\dagger}-a^{\dagger}a=1$ can be considered formally, in a suitable space of functions $f$, by letting $a=\frac{d}{dx}$ and $a^{\dagger}=x$ (the operator of multiplication by $x$). From \eqref{3}, we note that
\begin{equation}
\Big(x\frac{d}{dx}\Big)_{n,\lambda}f(x)=\sum_{k=1}^{n}S_{2,\lambda}(n,k)x^{k}\Big(\frac{d}{dx}\Big)^{k}f(x),\label{19}
\end{equation}
where $n$ is a positive integer. \par
Now, we consider the normal ordering of a degenerate integral power of the number operator $a^{\dagger}a$ in terms of the boson operators $a, a^{\dagger}$.
In view of \eqref{19}, the normal ordering of the degenerate $k$th power of the number operator $a^{\dagger}a$ in terms of boson operators $a, a^{\dagger}$ can be written in the form
\begin{equation}
(a^{\dagger}a)_{k,\lambda}=\sum_{l=0}^{k}S_{2,\lambda}(k,l)(a^{\dagger})^{l}a^{l},\quad (k\in\mathbb{N}).\label{20}
\end{equation}
From \eqref{11} and \eqref{20}, we note that
\begin{equation}
(a^{\dagger}a)_{k,\lambda}|m\rangle=(a^{\dagger}a)(a^{\dagger}a-\lambda)\cdots(a^{\dagger}a-(k-1)\lambda|m\rangle=(m)_{k,\lambda}|m\rangle,\label{21}	
\end{equation}
and
 \begin{equation}
(a^{\dagger}a)_{k,\lambda}|m\rangle=\sum_{l=0}^{k}S_{2,\lambda}(k,l)(a^{\dagger})^{l}a^{l}|m\rangle=\sum_{l=0}^{k}S_{2,\lambda}(k,l)(m)_{l}|m\rangle.\label{22}
 \end{equation}
Thus, by \eqref{21} and \eqref{22}, we get
\begin{equation}
(m)_{k,\lambda}=\sum_{l=0}^{k}S_{2,\lambda}(k,l)(m)_{l},\quad (k\ge 1).\label{23}	
\end{equation}
This is the classical expression for the degenerate $k$th power of $m$ in terms of the falling factorials $(m)_{l}$. This shows that \eqref{4} holds for all nonnegative integers $x=m=0,1,2,\dots$, which in turn implies \eqref{4} itself holds true.
On the other hand, by \eqref{11} and \eqref{20}, we get
\begin{align}
\langle z|(a^{\dagger}a)_{k,\lambda}|z\rangle &=\sum_{l=0}^{k}S_{2,\lambda}(k,l)\langle z|(a^{\dagger})^{l}a^{l}|z\rangle
=\sum_{l=0}^{k}S_{2,\lambda}(k,l)(\overline{z})^{l}z^{l}\langle z|z\rangle	\label{24} \\
&=\sum_{l=0}^{k}S_{2,\lambda}(k,l)|z|^{2l}=\phi_{k,\lambda}(|z|^{2}) \nonumber.
\end{align}
Let $f(t)=\langle z|e_{\lambda}^{a^{\dagger}a}(t)|z\rangle$. Then, by \eqref{24}, we get
\begin{align}
f(t)&=\langle z|e_{\lambda}^{a^{\dagger}a}(t)|z\rangle
=\sum_{k=0}^{\infty}\frac{t^{k}}{k!}\langle z|(a^{\dagger}a)_{k,\lambda}|z\rangle \label{25}\\
&=\sum_{k=0}^{\infty}\frac{t^{k}}{k!}\sum_{l=0}^{k}S_{2,\lambda}(k,l)|z|^{2l}=\sum_{k=0}^{\infty}\phi_{k,\lambda}(|z|^{2})\frac{t^{k}}{k!}.\nonumber
\end{align}
Indeed, the equation \eqref{25} says that $f(t)$ is the generating function of the degenerate Bell polynomials. To obtain an explicit expression for $f(t)$, we differentiate $f(t)$ with respect to $t$. It is not difficult to show that
\begin{equation}
a^{\dagger}a\,e_{\lambda}^{a^{\dagger}a-\lambda}(t)=e_{\lambda}^{a^{\dagger}a-\lambda}(t)a^{\dagger}a=a^{\dagger}e_{\lambda}^{aa^{\dagger}-\lambda}(t)a=a^{\dagger}e_{\lambda}^{a^{\dagger}a+1-\lambda}(t)a.\label{26}	
\end{equation}
From \eqref{25} and \eqref{26}, we note that
\begin{align}
\frac{\partial f(t)}{\partial t}&=\frac{\partial}{\partial t}\langle z|e_{\lambda}^{a^{\dagger}a}(t)|z\rangle=\langle z|a^{\dagger}a\, e_{\lambda}^{a^{\dagger}a-\lambda}(t)|z\rangle \label{27} \\
&=\langle z|a^{\dagger}e_{\lambda}^{a^{\dagger}a+1-\lambda}(t)a|z\rangle=e_{\lambda}^{1-\lambda}(t)\langle z|a^{\dagger}e_{\lambda}^{a^{\dagger}a}(t)a|z\rangle\nonumber\\
&=e_{\lambda}^{1-\lambda}(t)\overline{z}  z
\langle z|e_{\lambda}^{a^{\dagger}a}(t)|z\rangle=e_{\lambda}^{1-\lambda}(t)|z|^{2}f(t).\nonumber
\end{align}
Thus, we have
\begin{equation}
\frac{\partial f(t)}{\partial t}=e_{\lambda}^{1-\lambda}(t)|z|^{2}f(t)\ \Longleftrightarrow\ \frac{f^{\prime}(t)}{f(t)}=e_{\lambda}^{1-\lambda}(t)|z|^{2},\quad \bigg(f^{\prime}(t)=\frac{d}{dt}f(t)\bigg).\label{28}
\end{equation}
Assume that $f(0)=1$, for the initial value. Then, by \eqref{28}, we get
\begin{equation}
\log f(t)=\int_{0}^{t}\frac{f^{\prime}(t)}{f(t)}dt=\int_{0}^{t}|z|^{2}e_{\lambda}^{1-\lambda}(t)dt=\big(e_{\lambda}(t)-1\big)|z|^{2}. \label{29}	
\end{equation}
The equation \eqref{29} can be rewritten as
\begin{equation}
f(t)=e^{|z|^{2}(e_{\lambda}(t)-1)}=\sum_{l=0}^{\infty}|z|^{2l}\frac{1}{l!}(e_{\lambda}(t)-1)^{l}.\label{30}
\end{equation}
From \eqref{25}, we note that
\begin{align}
f(t)=\sum_{k=0}^{\infty}\frac{t^{k}}{k!}\bigg(\sum_{l=0}^{k}S_{2,\lambda}(k,l)|z|^{2l}\bigg)=\sum_{l=0}^{\infty}\bigg(\sum_{k=l}^{\infty}S_{2,\lambda}(k,l)\frac{t^{k}}{k!}\bigg)|z|^{2l}.\label{31}
\end{align}
Therefore, by \eqref{30} and \eqref{31}, we get
\begin{equation*}
\frac{1}{l!}\big(e_{\lambda}(t)-1\big)^{l}=\sum_{k=l}^{\infty}S_{2,\lambda}(k,l)\frac{t^{k}}{k!},
\end{equation*}
which agrees with \eqref{6}.
From \eqref{24}, we have
\begin{equation}
\langle z|(a^{\dagger}a)_{k,\lambda}|z\rangle=\sum_{l=0}^{k}S_{2,\lambda}(k,l)|z|^{2l}=\phi_{k,\lambda}(|z|^{2}).\label{32}	
\end{equation}
Setting $|z|=1$, we obtain
\begin{displaymath}
\langle z|(a^{\dagger}a)_{k,\lambda}|z\rangle=\phi_{k,\lambda},\quad (k\ge 1).
\end{displaymath}
Differentiating \eqref{25} with respect to $t$, we obtain
\begin{equation}
\frac{\partial f(t)}{\partial t}=\sum_{k=1}^{\infty}\frac{t^{k-1}}{(k-1)!}\phi_{k,\lambda}(|z|^{2})=\sum_{k=0}^{\infty}\phi_{k+1,\lambda}(|z|^{2})\frac{t^{k}}{k!}.\label{33}	
\end{equation}
On the other hand, by \eqref{28}, we get
\begin{align}
\frac{\partial f(t)}{\partial t}&=e_{\lambda}^{1-\lambda}(t)|z|^{2}f(t)=e_{\lambda}^{1-\lambda}(t)|z|^{2}\sum_{k=0}^{\infty}\frac{t^{k}}{k!}\phi_{k,\lambda}(|z|^{2})\label{34} \\
&=|z|^{2}\sum_{k=0}^{\infty}\bigg(\sum_{l=0}^{k}\binom{k}{l}(1-\lambda)_{k-l,\lambda}\phi_{l,\lambda}(|z|^{2})\bigg)\frac{t^{k}}{k!}. \nonumber	
\end{align}
Thus, by \eqref{33} and \eqref{34}, we get
\begin{equation*}
\phi_{k+1,\lambda}(|z|^{2})=|z|^{2}\sum_{l=0}^{k}\binom{k}{l}(1-\lambda)_{k-l,\lambda}\phi_{l,\lambda}(|z|^{2}).
\end{equation*}
In particular, for $|z|=1$, we have
\begin{equation*}
\phi_{k+1,\lambda}=\sum_{l=0}^{k}\binom{k}{l}(1-\lambda)_{k-l,\lambda}\phi_{l,\lambda}.
\end{equation*}
Evaluating the left hand side of \eqref{32} by using the representation of the coherent state in terms of the number state in \eqref{18}, we have
\begin{align}
\langle z|(a^{\dagger}a)_{k,\lambda}|z\rangle &=e^{-\frac{|z|^{2}}{2}}\cdot e^{-\frac{|z|^{2}}{2}}	\sum_{m,n=0}^{\infty}\frac{\bar{z}^{m}z^{n}}{\sqrt{m!}\sqrt{n!}}(n)_{k,\lambda}\langle m|n\rangle\label{35}\\
&=e^{-|z|^{2}}\sum_{n=0}^{\infty}\frac{|z|^{2n}}{n!}(n)_{k,\lambda}.	\nonumber
\end{align}
Thus, by \eqref{32} and \eqref{35}, we get
\begin{align}
\phi_{k,\lambda}(|z|^{2})&=\sum_{l=0}^{k}|z|^{2l}S_{2,\lambda}(k,l)
=e^{-|z|^{2}}\sum_{n=0}^{\infty}\frac{|z|^{2n}}{n!}(n)_{k,\lambda} \label{36} \\
&=e^{-|z|^{2}}\sum_{n=1}^{\infty}\frac{|z|^{2n}}{(n-1)!}(n-\lambda)_{k-1,\lambda},\quad (k\in\mathbb{N}). \nonumber
\end{align}
In particular, by letting $|z|=1$, we get
\begin{align}
\phi_{k,\lambda}=\frac{1}{e}\sum_{n=0}^{\infty}\frac{1}{n!}(n)_{k,\lambda}
=\frac{1}{e}\sum_{n=1}^{\infty}\frac{1}{(n-1)!}(n-\lambda)_{k-1,\lambda},\quad (k\in\mathbb{N}).\label{37}	
\end{align}
This is a Dobinski-like formula for the degenerate Bell numbers.

\section{conclusion}

Intensive studies have been done for degenerate versions of quite a few special polynomials and numbers by using such tools as combinatorial methods, generating functions, mathematical physics, umbral calculus techniques, $p$-adic analysis, differential equations, special functions, probability theory and analytic number theory.\par
As a degenerate version of the well known normal ordering of an integral power of the number operator, we considered the normal ordering of a degenerate integral power of the number operator in terms of boson operators. By using this normal ordering
we derived two equations defining the degenerate Stirling numbers of the second kind and a Dobinski-like formula for the degenerate Bell numbers. \par
It is one of our future projects to continue to explore various degenerate versions of many special polynomials and numbers by using aforementioned tools.

\vspace{0.5cm}

{ \bf  Disclosure statement}\

No potential conflict of interest was reported by the authors.

\vspace{0.5cm}

{\bf Funding}

	This work was supported by the Basic Science Research Program, the National Research Foundation of Korea, (NRF-2021R1F1A1050151).

\end{document}